\documentclass{ifacconf}

\usepackage{natbib} 
\usepackage[centertags]{amsmath}
\usepackage{graphics}
\usepackage{amscd}
\usepackage{amsfonts}
\usepackage{amssymb}
\usepackage{helvet}
\usepackage{rotating}
\usepackage{epsfig}
\usepackage[latin1]{inputenc}
\usepackage{longtable}
\usepackage{natbib}
\usepackage[english,german,american]{babel}
\usepackage{relsize}
\usepackage[european]{circuitikz}
\usepackage[chapter]{algorithm} % 
\usepackage{algorithmic}
\newtheorem{theorem}{Theorem}[section]

\newtheorem{proposition}[theorem]{Proposition}
\newtheorem{definition}[theorem]{Definition}
\newtheorem{assumption}[theorem]{Assumption}

\newtheorem{remark}[theorem]{Remark}
\newcommand{\R}{\mathbb{R}}
\newcommand{\N}{\mathbb{N}}

\newcommand{\cK}{\mathcal{K}}
\newcommand{\cKL}{\mathcal{KL}}

\newcommand{\state}{x}
\newcommand{\statestar}{\state^\star}
\newcommand{\statedisturbed}{\overline{\state}}
\newcommand{\statezero}{\state_0}
\newcommand{\statezerodisturbed}{\overline{\state}_0}
\newcommand{\statemax}{\Delta \state}
\newcommand{\control}{u}
\newcommand{\controlstar}{\control^\star}
\newcommand{\controldisturbed}{\overline{\control}}

\newcommand{\disturbancezero}{w}
\newcommand{\disturbancezerostar}{\disturbancezero^\star}
\newcommand{\disturbance}{\overline{\disturbancezero}}
\newcommand{\disturbancemax}{\Delta \disturbancezero}
\newcommand{\feedback}{\mu}
\newcommand{\feedbackdisturbed}{\overline{\feedback}}
\newcommand{\feedbackN}{\feedback_N}
\newcommand{\feedbackdisturbedN}{\feedbackdisturbed_N}

\newcommand{\distance}{d}
\newcommand{\dynamic}{f}
\newcommand{\stagecost}{\ell}
\newcommand{\costfunction}{J}
\newcommand{\costfunctionN}{\costfunction_N}
\newcommand{\costfunctioninfty}{\costfunction_\infty}
\newcommand{\valuefunction}{V}
\newcommand{\valuefunctionN}{\valuefunction_N}

\newcommand{\valuefunctioninfty}{\valuefunction_\infty}

\newcommand{\valuefunctioninftyfeedbackdisturbedN}{\valuefunctioninfty^{\feedbackdisturbedN}}
\newcommand{\closedloopindex}{n}
\newcommand{\openloopindex}{k}
\newcommand{\stateset}{X}
\newcommand{\controlset}{U}
\newcommand{\disturbanceset}{W}
\newcommand{\statesetconstrained}{\mathbb{\stateset}}
\newcommand{\controlsetconstrained}{\mathbb{\controlset}}
\newcommand{\openloopset}{\mathcal{I}}
\newcommand{\admissiblecontrolset}{\mathcal{\controlset}^\openloopset}
\newcommand{\updateset}{A}
\newcommand{\unknownset}{\mathcal{L}}

\begin{document}

% ------------------------------------------------------------------------
\begin{frontmatter}
\title{Performance bounds for NMPC combined with Sensitivity Updates}
\author{J\"{u}rgen Pannek, Matthias Gerdts}
\address{Faculty of Aeronautics, University of the Federal Armed Forces Munich, Germany; juergen.pannek@unibw.de, matthias.gerdts@unibw.de}

% ------------------------------------------------------------------------
\begin{abstract}
: In this paper we present a stability proof of model predictive control without stabilizing terminal constraints of cost which are subject to unknown but measurable disturbances. To this end, a relaxed Lyapunov argument on the nominal system and Lipschitz conditions on the open loop change of the optimal value function and the stage costs are employed. Based on the special case of sensitivity analysis, we show that Lipschitz assumptions are satisfied if a sensitivity update can be performed along the closed loop solution. To illustrate our approach we present a halfcar example and show performance improvement of the updated solution with respect to comfort and handling properties.
\end{abstract}
% ------------------------------------------------------------------------

\end{frontmatter}
% ------------------------------------------------------------------------

% ------------------------------------------------------------------------
\section{Introduction}
% ------------------------------------------------------------------------
Due to its simple structure, model predictive control (MPC) has become a well--established method for (sub)op\-ti\-mal control of linear and nonlinear systems, see, e.g., \cite{CamachoBordons2004} and \cite{RawlingsMayne2009, GruenePannek2011}. By means of this method an approximated closed--loop solution of an infinite horizon optimal control problem is computed in the following way: in each sampling interval, based on a measurement of the current state, a finite horizon optimal control problem is solved and the first element (or sometimes also more) of the resulting optimal control sequence is applied as input for the next sampling interval(s). This procedure is then repeated iteratively.

Unfortunately, stability and optimality may be lost due to the trunctation of the infinite horizon. In order to ensure stability of the resulting closed loop, one usually imposes terminal point constraints as shown in \cite{KeerthiGilbert1988} and \cite{Alamir2006} or Lyapunov type terminal costs and terminal regions, see \cite{ChenAllgoewer1998} and \cite{MayneRawlingsRaoScokaert2000}. A third approach uses a relaxed Lyapunov condition presented in \cite{GrueneRantzer2008} which can be shown to hold if the system is controllable in terms of the stage costs, cf. \cite{Gruene2009}. Additionally, this method allows for computing an estimate on the degree of suboptimality with respect to the infinite horizon controller, see also \cite{ShammaXiong1997} and \cite{NevisticPrimbs1997} for earlier works on this topic.

Here, we follow the third approach but extend its applicability for the case of parametric control systems and subsequent disturbance rejection updates. In particular, we impose an abstract update of the control law with respect to the measured disturbance or parameter. Then, if Lipschitz conditions on the open loop change on both the optimal value and stage cost function hold for such an update, we can utilize the relaxed Lyapunov condition for the nominal open loop solutions to show stability of the resulting closed loop. 

Thereafter, we focus on the update law for optimal controls via sensitivities, see also \cite{ZavalaBiegler2009}. Using results from sensitivity analysis, we can show that if a sensitivity update can be performed along the closed loop, then the Lipschitz assumptions hold and stability of the disturbed system can be guaranteed. Note that due to abstracting from the form of the update, the proposed stability proof is not limited to sensitivities but may also applicable for other update methods such as realtime iterations or hierarchical MPC, see \cite{DiehlBockSchloeder2005} and \cite{BockDiehlKostinaSchloeder2007} respectively.

The paper is organized as follows. In Section \ref{Section:Preliminaries} the problem formulation and the concept of practical stability are defined. In the subsequent Section \ref{Section:Stability} the stability proof of MPC subject to disturbances without stabilizing terminal constraints and costs is given. In Section \ref{Section:Sensitivity theory}, the additional assumptions required within the stability proof are shown to be satisfied in case of sensitivity based updates of the control law. Moreover, we present simulation results for a halfcar subject to disturbed measurements of both the state and the sensor inputs in Section \ref{Section:Numerical Results}. To conclude our paper, we draw some conclusions and give an outlook on future research.

% ------------------------------------------------------------------------
\section{Setup and Preliminaries}
\label{Section:Preliminaries}
% ------------------------------------------------------------------------

The types of control systems we consider within this work are given by the dynamics
\begin{align}
	\label{Preliminaries:eq:control system}
	\state(\closedloopindex+1) = \dynamic(\state(\closedloopindex),\control(\closedloopindex), \disturbancezero(\closedloopindex))
\end{align}
where $\state$ denotes the state of the system, $\control$ the external control and $\disturbancezero$ an disturbance which can be measured. These variables are elements of respective metric spaces $(\stateset, \distance_{\stateset})$, $(\controlset, \distance_{\controlset})$ and $(\disturbanceset, \distance_{\disturbanceset})$ which represent the state, control and disturbance space. Therefore, our results are also applicable to discrete time dynamics induced by a sampled finite or infinite dimensional system. For ease of notation, we introduce the abbreviation $\| x \|_{y} = \distance_i(x, y)$ for $\distance_i \in \{ \distance_\stateset, \distance_\controlset, \distance_\disturbanceset \}$. Both, state and control, are constrained to be elements of subsets $\statesetconstrained \subseteq \stateset$ and $\controlsetconstrained \subseteq \controlset$. We denote the \textit{undisturbed} or \textit{nominal} state trajectory corresponding to an initial state $\statezero \in \stateset$, a control sequence $\control = \left( \control(\openloopindex) \right)_{\openloopindex \in \openloopset}$ and a nominal disturbance sequence $\disturbancezero = \left( \disturbancezero(\openloopindex) \right)_{\openloopindex \in \openloopset}$, $\openloopset := \{ 0, 1, \ldots, N - 1\}$ with $N \in \N$ or $\openloopset := \N_0$ by $\state(\cdot) = \state_{\control, \disturbancezero}(\cdot; \statezero)$. Similarly, the \textit{disturbed} state trajectory which is subject to a disturbed initial value $\statezerodisturbed$, a control sequence $\controldisturbed = \left( \controldisturbed(\openloopindex) \right)_{\openloopindex \in \openloopset}$ the disturbance sequence $\disturbance = \left( \disturbance(\openloopindex) \right)_{\openloopindex \in \openloopset}$ is denoted by $\statedisturbed(\cdot) = \statedisturbed_{\controldisturbed, \disturbance}(\cdot; \statezerodisturbed)$. Since in the presence of constraints not all control sequences $\control$ are admissible, we introduce $\admissiblecontrolset(\statezero, \disturbancezero)$ as the \textit{set of all admissible control sequences} $\control = \left( \control(\openloopindex) \right)_{\openloopindex \in \openloopset}$ which for a fixed disturbance sequence $\disturbancezero$ satisfy $\dynamic (\state(\openloopindex), \control(\openloopindex), \disturbancezero(\openloopindex)) \in \statesetconstrained$ and $\control(\openloopindex) \in \controlsetconstrained$ for $\openloopindex \in \openloopset$. 

Our aim in this work is to stabilize system \eqref{Preliminaries:eq:control system} at a controlled nominal equilibrium, that is a point $\statestar$ such that there exists a control $\controlstar$ and a nominal disturbance $\disturbancezerostar$ satisfying $\dynamic(\statestar, \controlstar, \disturbancezerostar) = \statestar$. To this end, we consider a two stage feedback design. In the first stage we want to compute a static state feedback law for a given nominal disturbance sequence $\disturbancezero = \left( \disturbancezero(\closedloopindex) \right)_{\closedloopindex \in \N_0}$ which minimizes the infinite horizon cost functional $\costfunctioninfty ( \state, \control, \disturbancezero ) = \sum_{\closedloopindex=0}^\infty \stagecost(\state(\closedloopindex), \control(\closedloopindex), \disturbancezero(\closedloopindex))$. In this context, the stage costs $\stagecost: \stateset \times \controlset \times \disturbanceset \to \R_0^+$ are continuous and satisfy $\stagecost(\statestar, \controlstar, \disturbancezerostar) = 0$ and $\stagecost(\state, \control, \disturbancezero) > 0$ for all $\control \in \controlset$ for each $\state \not = \statestar$ and each $\disturbancezero \not = \disturbancezerostar$. In order to avoid solving the discrete time equivalent to a Hamilton--Jacobi--Bellman equation which is computationally intractable in most cases, we use a model predictive control (MPC) approach to approximate the desired solution. The resulting cost functional is given by
\begin{align}
	\label{Preliminaries:eq:cost functional}
	\costfunctionN(\state, \control, \disturbancezero) := \sum_{\openloopindex=0}^N \stagecost(\state(\openloopindex; \state), \control(\openloopindex), \disturbancezero(\openloopindex))
\end{align}
where $N \in \N_{\geq 2}$ denotes the length of the prediction horizon, i.e. the prediction horizon is truncated and, thus, finite. Consequently, the obtained control sequence $\controlstar(\cdot, \state, \disturbancezero)$ itself is also finite. \\
In order to retrieve an infinite sequence, one usually implements only the first part of this sequence $\feedbackN(\state, \disturbancezero) := \controlstar(0, \state, \disturbancezero)$, then the prediction horizon is shifted forward in time and the procedure is iterated ad infimum, cf., e.g., \cite{CamachoBordons2004, Alamir2006} or \cite{RawlingsMayne2009}.
\begin{remark}
	Note that even the evaluation of both $\costfunctioninfty$ and $\costfunctionN$ require knowledge of a sequence of nominal future disturbances. Such values may be obtained by forward measurements, e.g. cameras in a car or thermometers at intake pipes, or extrapolation methods.
\end{remark}
Here, we want to update the feedback law depending on intermediate disturbances $\disturbance$ and newly obtained state estimates $\statedisturbed$. Hence,  while shifting and iterating the problem remains unchanged, a modified control $\feedbackdisturbedN(\statedisturbed, \disturbance)$ is implemented instead of $\feedbackN(\state, \disturbancezero)$ and we assume $\feedbackdisturbedN(\statedisturbed, \disturbance)$ to be instantly computable.

For simplicity of exposition, we assume that a minimizer $\controlstar(\cdot, \state, \disturbancezero)$ of \eqref{Preliminaries:eq:cost functional} exists for each $\state \in \statesetconstrained$ and $N \in \N$. Particularly, this includes the assumption that a feasible solution exists for each $\state \in \statesetconstrained$ and each $\disturbancezero \in \disturbanceset$. For methods on avoiding the feasibility assumption in $\state$ we refer to \cite{PrimbsNevistic2000} or \cite{GruenePannek2011}. Using the existence of a minimizer $\controlstar(\cdot, \state, \disturbancezero) \in \admissiblecontrolset(\state, \disturbancezero)$, we obtain the following equality for the optimal value function defined on a finite horizon
\begin{align}
	\label{Preliminaries:eq:value function}
	\valuefunctionN(\state, \disturbancezero) = \inf_{\control \in \admissiblecontrolset(\state, \disturbancezero)} \costfunctionN(\state, \control, \disturbancezero).% = \costfunctionN(\state,\controlstar(\cdot, \disturbancezero), \disturbancezero).
\end{align}
In order to compute a performance or suboptimality index of the updated MPC feedback $\feedbackdisturbedN(\statedisturbed, \disturbance)$, we denote the closed loop trajectory by
\begin{align*}
	\statedisturbed(\closedloopindex+1) = \dynamic(\statedisturbed(\closedloopindex), \feedbackdisturbedN(\statedisturbed(\closedloopindex), \disturbance(\closedloopindex)), \disturbance(\closedloopindex))
\end{align*}
which gives rise to the closed loop costs
\begin{align}
	\valuefunctioninftyfeedbackdisturbedN(\statedisturbed, \disturbance)  := \sum_{\closedloopindex=0}^{\infty} \stagecost(\statedisturbed(\closedloopindex), \feedbackdisturbedN(\statedisturbed(\closedloopindex), \disturbance(\closedloopindex)), \disturbance(\closedloopindex)).
\end{align}
Regarding the disturbances, we assume the following: 
\begin{assumption}\label{Preliminaries:ass:bounded disturbance}
	The disturbance is bounded by $\disturbancemax > 0$ and the maximal state estimate deviation is bounded by $\statemax > 0$, i.e. $\| \disturbance \|_{\disturbancezero} \leq \disturbancemax$ and $\| \statedisturbed \|_{\state} \leq \statemax$. Additionally, the nominal disturbance change is bounded by $\disturbancemax$, that is $\| \disturbance(n+1) \|_{\disturbance(n)} \leq \disturbancemax$.
\end{assumption}
Due to the presence of unknown disturbances, asymptotic stability cannot be expected. Therefore, we consider the concept of $P$--practical asymptotic stability. Similar to input--to--state stability (ISS), the convergence of the closed loop solution is characterized by comparison functions. Here, we call a continuous function $\rho: \R_{\geq 0} \rightarrow \R_{\geq 0}$ a class $\cK_\infty$-function if it satisfies $\rho(0) = 0$, is strictly increasing and unbounded. A continuous function $\beta: \R_{\geq 0} \times \R_{\geq 0} \rightarrow \R_{\geq 0}$ is said to be of class $\cKL$ if for each $r > 0$ the limit $\lim_{t \rightarrow \infty} \beta(r, t) = 0$ holds and for each $t \geq 0$ the condition $\beta(\cdot, t) \in \cK_\infty$ is satisfied. 
\begin{definition}\label{Preliminaries:def:p practical}
	Let $\updateset \subset \stateset$ be a forward invariant set with respect to all possible disturbances satisfying Assumption \ref{Preliminaries:ass:bounded disturbance} and let $P \subset \updateset$. Then $\statestar \in P$ is $P$--practically asymptotically stable on $\updateset$ if there exists $\beta \in \cKL$ such that
	\begin{align}
		\label{Preliminaries:eq:closed loop}
		\| \statedisturbed(\closedloopindex) \|_{\statestar} \leq \beta( \| \statezerodisturbed \|_{\statestar}, \closedloopindex)
	\end{align}
	holds for all $\statezerodisturbed \in \updateset$ and all $n \in \N_0$ with $\statedisturbed(\closedloopindex) \not \in P$.
\end{definition}
Note that the ISS property can be shown for $P$--practical asymptotically stable systems by a suitable choice of the comparison function $\gamma \in \cK$, cf. Chapter 8.5 in \cite{GruenePannek2011}. Typically, the ISS property is shown via ISS Lyapunov functions, see, e.g, \cite{JiangWang2001} or \cite{MagniScattolini2007}.
Different from that, $P$--practical asymptotic stability can be concluded if there exists a suitable ``truncated'' Lyapunov function as shown in Theorem 2.20 in \cite{GruenePannek2011}.
\begin{theorem}\label{Preliminaries:thm:p practical}
	Suppose  $\updateset \subset \stateset$ is a forward invariant set with respect to all possible disturbances satisfying Assumption \ref{Preliminaries:ass:bounded disturbance}, $P \subset \updateset$ and $\statestar \in P$. If there exist $\cK$--functions $\alpha_1$, $\alpha_2$, $\alpha_3$ and a Lyapunov function $\valuefunction$ on $S = \updateset \setminus P$ satisfying
	\begin{align*}
		\alpha_1( \| \statedisturbed \|_{\statestar} + \| \disturbance \|_{\disturbancezerostar} ) & \leq \valuefunction(\statedisturbed, \disturbance) \leq  \alpha_2( \| \statedisturbed \|_{\statestar} + \| \disturbance \|_{\disturbancezerostar} ) \\
		\valuefunction(\statedisturbed, \disturbance) & \geq \valuefunction(\dynamic(\statedisturbed, \feedbackdisturbedN(\statedisturbed, \disturbance), \disturbance), \disturbance) \\
		& \qquad - \alpha_3( \| \statedisturbed \|_{\statestar} )
	\end{align*}
	then $\statestar$ is $P$--practically asymptotically stable on $\updateset$.
\end{theorem}

% ------------------------------------------------------------------------
\section{Stability}
\label{Section:Stability}
% ------------------------------------------------------------------------

In the literature, one usually uses terminal constraints or Lyapunov type terminal costs to guarantee the ISS property of the disturbed closed loop, cf., e.g., \cite{DiehlBockSchloeder2005, BockDiehlKostinaSchloeder2007, ZavalaBiegler2009}. Here, we consider the plain MPC formulation without these modifications. Instead, we suppose a relaxed Lyapunov condition to hold for the nominal case which additionally reveals a performance index of the nominal closed loop, cf. \cite{LincolnRantzer2006} and \cite{GrueneRantzer2008}. Note that this condition is always satisfied if $N$ is sufficiently large, see \cite{AlamirBornard1995}, \cite{JadbabaieHauser2005} or \cite{GrimmMessinaTeel2005}.

In order to prove a similar result in the disturbed case, one could modify the stage cost $\stagecost$ to be positive definite with respect to a robustly stabilizable forward invariant neighbourhood of $\statestar$. Since the computation of this neighbourhood may be impossible, we choose the stage cost $\stagecost$ to be positive definite with respect to $\statestar$ only, that is ignoring the effects of disturbances on stabilizability. As $\stagecost$ is typically much smaller close to $\statestar$ than far away from the desired steady state, we may still expect the closed loop to converge to a neighbourhood of $\statestar$, i.e. $P$--practical stability of the closed loop. Similar to results shown in \cite{GrueneRantzer2008} we additionally obtain a bound for the degree of suboptimality.

In order to show such a performance results, we assume the following:
\begin{assumption}\label{Stability:ass:sensitivity}
	There exists a set $\updateset \subset \statesetconstrained$ containing $\statestar$ such that
	\begin{align*}
		%& \distance_{\controlset} ( \feedbackdisturbedN(\statedisturbed, \disturbance), \feedbackdisturbedN(\state, \disturbancezero)) ) \leq L_\control \left( \| \statedisturbed \|_{\state} + \| \disturbance \|_{\disturbancezero} \right) \\ \displaybreak[0]
		& | \stagecost(\statedisturbed, \feedbackdisturbedN(\statedisturbed, \disturbance), \disturbance) - \stagecost(\state, \feedbackN(\state, \disturbancezero), \disturbancezero) | \\
		& \quad \leq L_\stagecost \left( \| \statedisturbed \|_{\state} + \| \disturbance \|_{\disturbancezero} \right) \\ \displaybreak[0]
		& | \costfunctionN( \statedisturbed, \feedbackdisturbedN(\statedisturbed, \disturbance), \disturbance) - \costfunctionN(\state, \feedbackN(\state, \disturbancezero), \disturbancezero) | \\
		& \quad \leq L_\costfunction \left( \| \statedisturbed \|_{\state} + \| \disturbance \|_{\disturbancezero} \right)
	\end{align*}
	hold with Lipschitz constants %$L_\control$, 
	$L_\stagecost$ and $L_\costfunction$ for all tupels $(\statedisturbed, \state, \disturbance, \disturbancezero)$ with $\statedisturbed, \state \in \updateset$, $\disturbance, \disturbancezero \in \disturbanceset$ satisfying Assumption \ref{Preliminaries:ass:bounded disturbance}.% Additionally, the dynamic $\dynamic$ is Lipschitz continuous with constant $L_\dynamic$ with respect to all variables.
\end{assumption}
Under these conditions, we can prove the following:
\begin{theorem}\label{Stability:thm:performance}
	Suppose Assumptions \ref{Preliminaries:ass:bounded disturbance}, \ref{Stability:ass:sensitivity} to hold and a given feedback $\feedbackN: \stateset \times \disturbanceset \to \controlset$ and a nonnegative function $\valuefunctionN: \stateset \times \disturbanceset \to \R_0^+$ to satisfy the relaxed Lyapunov inequality
	\begin{align}
		\label{Stability:thm:performance:eq1}
		\valuefunctionN(\state, \disturbancezero) \geq \valuefunctionN(\dynamic(\state,\feedbackN(\state, \disturbancezero), \disturbancezero), \disturbancezero) + \alpha \ell(\state, \feedbackN(\state, \disturbancezero), \disturbancezero )
	\end{align}
	for some $\alpha \in (0, 1)$ and all $\state \in \updateset$ and all $\disturbancezero \in \disturbanceset$. Suppose $\varepsilon \geq L_\stagecost ( \statemax + \disturbancemax) + ( 2 L_\costfunction \statemax + 3 L_\costfunction \disturbancemax) / \alpha$ and let $\unknownset$ denote the minimal set which contains $\statestar$, is forward invariant with respect to all possible disturbances satisfying Assumption \ref{Preliminaries:ass:bounded disturbance} and $\valuefunctionN(\state, \disturbancezero) \geq \valuefunctionN(\dynamic(\state,\feedbackN(\state, \disturbancezero), \disturbancezero) + \alpha \varepsilon$ holds for all $\state \in \updateset \setminus \unknownset$ and all $\disturbancezero \in \disturbanceset$. Furthermore define the modified costs
	\begin{align}
		\label{Stability:thm:performance:eq2}
		\overline{\stagecost}(\state, \control, \disturbancezero) := \begin{cases}
			\max\left\{\stagecost( \state, \control, \disturbancezero) - \varepsilon, 0\right\} & \state \in \updateset \setminus \unknownset \\
			0 & \state \in \unknownset
		\end{cases}
	\end{align}
	and $\sigma :=\inf\{\valuefunctionN(\dynamic(\state, \feedbackN(\state, \disturbancezero), \disturbancezero), \disturbancezero) \mid \state \in \updateset \setminus \unknownset, \disturbancezero \in \disturbanceset \}$. Then for the modified closed loop cost
	\begin{align*}
		\overline{\valuefunctioninftyfeedbackdisturbedN}(\statedisturbed, \disturbance)   := \sum_{\closedloopindex=0}^{\infty} \overline{\stagecost}(\statedisturbed(\closedloopindex), \feedbackdisturbedN(\statedisturbed(\closedloopindex), \disturbance(\closedloopindex)), \disturbance(\closedloopindex))
	\end{align*}
	we have
	\begin{align}
		\label{Stability:thm:performance:eq3}
		\alpha \overline{\valuefunctioninftyfeedbackdisturbedN}(\statedisturbed, \disturbance) \leq \valuefunctionN(\statedisturbed, \disturbance) - \sigma \leq \valuefunctioninfty(\statedisturbed, \disturbance) - \sigma
	\end{align}
	for all $\statedisturbed \in \updateset$.
\end{theorem}
\textbf{Proof:} Consider $\statezero \in \updateset$. Let $\closedloopindex_0 \in \N_0$ be minimal with $\statedisturbed(\closedloopindex_0 + 1) \in \unknownset$ and set $\closedloopindex_0 := \infty$ if this case does not occur.\\
	Reformulating \eqref{Stability:thm:performance:eq1} we obtain
	\begin{align*}
				\alpha \stagecost(\state, \feedbackN(\state, \disturbancezero), \disturbancezero) \leq \valuefunctionN(\state, \disturbancezero) - \valuefunctionN(\dynamic(\state,\feedbackN(\state, \disturbancezero), \disturbancezero), \disturbancezero).
	\end{align*}
	Now we incorporate the effects of disturbances and sensitivity updates of the control using Assumption \ref{Stability:ass:sensitivity} which gives us
	\begin{align*}
		& \alpha \stagecost(\statedisturbed(\closedloopindex), \feedbackdisturbedN(\statedisturbed(\closedloopindex), \disturbance(\closedloopindex) ), \disturbance(\closedloopindex)) \\ \displaybreak[0]
		& \leq \alpha \stagecost(\state(\closedloopindex), \feedbackN(\state(\closedloopindex), \disturbancezero(\closedloopindex)), \disturbancezero(\closedloopindex)) \\
		& \qquad + \alpha L_\stagecost \left( \| \statedisturbed(\closedloopindex) \|_{\state(\closedloopindex)} + \| \disturbance(\closedloopindex) \|_{\disturbancezero(\closedloopindex)} \right) \\ \displaybreak[0]
		& \leq \valuefunctionN(\state(\closedloopindex), \disturbancezero(\closedloopindex)) +  \alpha L_\stagecost \left( \| \statedisturbed(\closedloopindex) \|_{\state(\closedloopindex)} + \| \disturbance(\closedloopindex) \|_{\disturbancezero(\closedloopindex)} \right) \\
		& \qquad - \valuefunctionN(\dynamic(\state(\closedloopindex),\feedbackN(\state(\closedloopindex), \disturbancezero(\closedloopindex)), \disturbancezero(\closedloopindex)), \disturbancezero(\closedloopindex)) \\ \displaybreak[0]
		& \leq \valuefunctionN(\state(\closedloopindex), \disturbancezero(\closedloopindex)) +  \alpha L_\stagecost \left( \| \statedisturbed(\closedloopindex) \|_{\state(\closedloopindex)} + \| \disturbance(\closedloopindex) \|_{\disturbancezero(\closedloopindex)} \right) \\
		& \qquad - \valuefunctionN(\dynamic(\state(\closedloopindex),\feedbackN(\state(\closedloopindex), \disturbancezero(\closedloopindex)), \disturbancezero(\closedloopindex)), \disturbancezero(\closedloopindex + 1)) \\ 
		& \qquad + L_\costfunction \| \disturbancezero(\closedloopindex+1) \|_{\disturbancezero(\closedloopindex)} \\ \displaybreak[0]
		& \leq \valuefunctionN(\state(\closedloopindex), \disturbancezero(\closedloopindex)) - \valuefunctionN(\state(\closedloopindex + 1), \disturbancezero(\closedloopindex + 1)) \\
		& \qquad + \alpha L_\stagecost \left( \| \statedisturbed(\closedloopindex) \|_{\state(\closedloopindex)} + \| \disturbance(\closedloopindex) \|_{\disturbancezero(\closedloopindex)} \right) \\
		& \qquad + L_\costfunction \| \disturbancezero(\closedloopindex+1) \|_{\disturbancezero(\closedloopindex)} \\ \displaybreak[0]
		& \leq \valuefunctionN(\statedisturbed(\closedloopindex), \disturbance(\closedloopindex)) - \valuefunctionN(\statedisturbed(\closedloopindex+1), \disturbance(\closedloopindex + 1)) \\
		& \qquad + ( \alpha L_\stagecost + L_\costfunction ) \left( \| \statedisturbed(\closedloopindex) \|_{\state(\closedloopindex)} + \| \disturbance(\closedloopindex) \|_{\disturbancezero(\closedloopindex)} \right) \\
		& \qquad + L_\costfunction \left( \| \statedisturbed(\closedloopindex+1) \|_{\state(\closedloopindex+1)} + \| \disturbance(\closedloopindex+1) \|_{\disturbancezero(\closedloopindex+1)} \right) \\
		& \qquad + L_\costfunction \| \disturbancezero(\closedloopindex+1) \|_{\disturbancezero(\closedloopindex)}
	\end{align*}
	Hence, using boundedness from Assumption \ref{Preliminaries:ass:bounded disturbance} reveals
	\begin{align*}
		& \alpha \stagecost(\statedisturbed(\closedloopindex), \feedbackdisturbedN(\statedisturbed(\closedloopindex), \disturbance(\closedloopindex)), \disturbance(\closedloopindex)) \\
		& \leq \valuefunctionN(\statedisturbed(\closedloopindex), \disturbance(\closedloopindex)) - \valuefunctionN(\statedisturbed(\closedloopindex+1), \disturbance(\closedloopindex + 1)) \\
		& \qquad + ( \alpha L_\stagecost + 2 L_\costfunction ) \statemax + ( \alpha L_\stagecost + 3 L_\costfunction ) \disturbancemax.
	\end{align*}
	Since we have that $\valuefunctionN(\state, \disturbancezero) \geq \valuefunctionN(\dynamic(\state,\feedbackN(\state, \disturbancezero), \disturbancezero), \disturbancezero) + ( \alpha L_\stagecost + 2 L_\costfunction ) \statemax + ( \alpha L_\stagecost + 3 L_\costfunction ) \disturbancemax$ holds for all $\state \in \updateset \setminus \unknownset$, we obtain $\valuefunctionN(\statedisturbed, \disturbance) \geq \valuefunctionN(\dynamic(\statedisturbed,\feedbackdisturbedN(\statedisturbed, \disturbance), \disturbance), \disturbance)$ for all $\state \in \updateset \setminus \unknownset$ and all $\disturbance$ satisfying Assumption \ref{Preliminaries:ass:bounded disturbance}. Using this fact, the bound on $\varepsilon$ and the definition of $\overline{\stagecost}$ in \eqref{Stability:thm:performance:eq2} we have
	\begin{align*}
		& \alpha \overline{\stagecost}(\statedisturbed(\closedloopindex), \feedbackdisturbedN(\statedisturbed(\closedloopindex), \disturbance(\closedloopindex)), \disturbance(\closedloopindex)) \\
		& = \max\left\{ \alpha \stagecost(\statedisturbed(\closedloopindex), \feedbackdisturbedN(\statedisturbed(\closedloopindex), \disturbance(\closedloopindex)), \disturbance(\closedloopindex)) - \alpha \varepsilon, 0\right\} \\
		& \leq \valuefunctionN(\statedisturbed(\closedloopindex), \disturbance(\closedloopindex)) - \valuefunctionN(\statedisturbed(\closedloopindex+1), \disturbance(\closedloopindex+1)).
	\end{align*}
	For $\closedloopindex \geq \closedloopindex_0 + 1$ the invariance of $\unknownset$ gives us $\statedisturbed(\closedloopindex) \in \unknownset$ and hence $\overline{\stagecost}(\statedisturbed(\closedloopindex), \feedbackdisturbedN(\statedisturbed(\closedloopindex), \disturbance(\closedloopindex)), \disturbance(\closedloopindex)) = 0$. Additionally, since $\sigma$ is the minimal cost after entry in $\unknownset$, we have $\valuefunctionN(\statedisturbed(\closedloopindex), \disturbance(\closedloopindex)) \geq \sigma$ for all $\closedloopindex \leq \closedloopindex_0 + 1$. Now we can sum the stage costs over $\closedloopindex$ and obtain
	\begin{align*}
		& \alpha \sum_{\closedloopindex = 0}^{K} \overline{\stagecost}(\statedisturbed(\closedloopindex), \feedbackdisturbedN(\statedisturbed(\closedloopindex), \disturbance(\closedloopindex)), \disturbance(\closedloopindex)) \\
		& = \alpha \sum_{\closedloopindex = 0}^{K_0} \overline{\stagecost}(\statedisturbed(\closedloopindex), \feedbackdisturbedN(\statedisturbed(\closedloopindex), \disturbance(\closedloopindex)), \disturbance(\closedloopindex)) \\
		& \leq \valuefunctionN(\statedisturbed(0), \disturbance(0)) - \valuefunctionN(\statedisturbed(K_0 + 1) \disturbance(K_0 + 1)) \\
		& \leq \valuefunction(\statedisturbed(0), \disturbance(0)) - \sigma.
	\end{align*}
	where $K_0 := \min \{ K, \closedloopindex_0 \}$. Using that $K \in \N$ was arbitrary we can conclude that $(\valuefunctionN(\statedisturbed(0), \disturbance(0)) - \sigma)/\alpha$ is an upper bound for $\overline{\valuefunctioninftyfeedbackdisturbedN}(\statedisturbed, \disturbance)$.\qed
	
The definition of $\updateset$ in Theorem \ref{Stability:thm:performance} is implicit, yet an approximation of $\updateset$ can be obtained via techniques presented in \cite{GrimmMessinaTeel2005} or \cite{GrueneRantzer2008}.

In addition to the previous performance estimate, $P$--practical asymptotic stability can be shown as follows:
\begin{theorem}
	Suppose the conditions of Theorem \ref{Stability:thm:performance} hold and additionally there exist $\cK$ functions $\alpha_1$, $\alpha_2$ such that
	\begin{align*}
		\alpha_1( \| \statedisturbed \|_{\statestar} + \| \disturbance \|_{\disturbancezerostar} ) & \leq \valuefunction(\statedisturbed, \disturbance) \leq  \alpha_2( \| \statedisturbed \|_{\statestar} + \| \disturbance \|_{\disturbancezerostar} )
	\end{align*}
	holds for all $\statedisturbed \in \updateset \setminus P$ with $P = \unknownset$. Then $\statestar$ is $P$--practically asymptotically stable on $\updateset$.
\end{theorem}
\textbf{Proof:} Follow directly from the definition of $\overline{\stagecost}$ in Theorem \ref{Stability:thm:performance} and the property $\valuefunction(\statedisturbed, \disturbance) \geq \valuefunctionN(\dynamic(\statedisturbed,\feedbackdisturbedN(\statedisturbed, \disturbance), \disturbance), \disturbance)$ shwon in the proof of Theorem \ref{Stability:thm:performance}. \qed

Note that the stability result holds for any update on the feedback law $\feedbackdisturbedN$ and is not specific to how this update is obtained.
\begin{remark}
	Condition \eqref{Stability:thm:performance:eq1} can be relaxed to
	\begin{align*}
		& \valuefunctionN(\state, \disturbancezero) - \valuefunctionN(\dynamic(\state,\feedbackN(\state, \disturbancezero), \disturbancezero), \disturbancezero) \\
		& \geq \min \{ \alpha \ell(\state, \feedbackN(\state, \disturbancezero), \disturbancezero ) - \overline{\varepsilon}, \ell(\state, \feedbackN(\state, \disturbancezero), \disturbancezero ) - \overline{\varepsilon} \}
	\end{align*}
	to additionally allow the nominal system to be $P$--practically stable only, cf. \cite{GrueneRantzer2008}. In this case, we have to modifiy the lower bound on $\varepsilon$ to $( ( L_\ell + 2 L_\costfunction ) (\disturbancemax + \statemax) + \overline{\varepsilon} ) / \alpha$ and require $\valuefunctionN(\state, \disturbancezero) \geq \valuefunctionN(\dynamic(\state,\feedbackN(\state, \disturbancezero), \disturbancezero) + 2 L_\costfunction (\disturbancemax + \statemax) + \overline{\varepsilon}$ to hold for all $\state \in \updateset \setminus \unknownset$.
\end{remark}
After stating stability for an abstract update law of $\feedbackdisturbedN$, our goal in the next section is to verify the required Assumption \ref{Stability:ass:sensitivity} for a particular updating strategy.

% ------------------------------------------------------------------------
\section{Sensitivity theory}
\label{Section:Sensitivity theory}
% ------------------------------------------------------------------------

Sensitivity analysis has been analysed extensively for the case of open loop optimal controls, see, e.g, \cite{GroetschelKrumkeRambau2001}, and has become a rather popular control method. In the MPC feedback context, \cite{ZavalaBiegler2009} analysed the impact of sensitivity updates on stability of the closed loop in the presence of stabilizing terminal constraints and Lyapunov type terminal costs. Here, we use an identical update
\begin{align}
	\label{Sensitivity theory:eq:update formula}
	\controldisturbed(\cdot, \statedisturbed, \disturbance) := \controlstar(\cdot, \state, \disturbancezero) + 
	\begin{pmatrix}
		\frac{\partial \controlstar}{\partial \state}(\cdot, \state, \disturbancezero) \\
		\frac{\partial \controlstar}{\partial \disturbancezero}(\cdot, \state, \disturbancezero)
	\end{pmatrix}^\top
	\begin{pmatrix}
		\statedisturbed(\cdot) - \state(\cdot) \\
		\disturbance(\cdot) - \disturbancezero(\cdot)
	\end{pmatrix}
\end{align}
and set $\feedbackdisturbedN(\statedisturbed, \disturbance) = \controldisturbed(0, \statedisturbed, \disturbance)$ but consider the plain MPC case described in Section \ref{Section:Preliminaries}. Since the solution of the underlying optimal control problem and the computation of the sensitivities require some computing time, such an approach is typically implemented in an advanced step setting, see also \cite{FindeisenAllgoewer2004}. The idea of the advanced step is to precompute an open loop control $\controlstar(\cdot, \state, \disturbancezero)$ and the sensitivities $\frac{\partial \controlstar}{\partial \state}$, $\frac{\partial \controlstar}{\partial \disturbancezero}$  for a future time instant while the current sampling period evolves. Then, once the time instant is reached at which the control is supposed to be implemented, newly obtained measurements of $\statedisturbed$ and $\disturbance$ are used to update the control according to \eqref{Sensitivity theory:eq:update formula}. An algorithmic implementation of this idea is shown below.
\begin{algorithm}\label{Sensitivity theory:alg:mpc}
	MPC Algorithm\\
	\begin{enumerate}
		\item Obtain measurement of $\statedisturbed(\closedloopindex)$ and $\disturbance(\closedloopindex)$
		\item Update control $\controldisturbed(\cdot, \statedisturbed(\closedloopindex), \disturbance(\closedloopindex)$ (e.g. via \eqref{Sensitivity theory:eq:update formula}) and apply control $\feedbackdisturbedN(\statedisturbed(\closedloopindex), \disturbance(\closedloopindex)) = \controldisturbed(0, \statedisturbed(\closedloopindex), \disturbance(\closedloopindex)$
		\item Predict $\state(\closedloopindex+1)$ using dynamic \eqref{Preliminaries:eq:control system} together with $\statedisturbed(\closedloopindex)$, $\feedbackdisturbedN(\statedisturbed(\closedloopindex), \disturbance(\closedloopindex))$ and $\disturbance(\closedloopindex)$
		\item Compute $\controlstar(\cdot, \state(\closedloopindex+1), \disturbance(\closedloopindex))$, set $\closedloopindex := \closedloopindex + 1$ and goto step (1)
	\end{enumerate}
\end{algorithm}

The mathematical foundations of such a control update are given in \cite{Fiacco1983}. For notational convenience, we adapted these results to the considered MPC case.
\begin{theorem}\label{Sensitivity theory:thm:sensitivity}
	Suppose that $\dynamic$ and $\stagecost$ are twice continuously differentiable in a neighbourhood of the nominal solution $\controlstar(\cdot, \state, \disturbancezero)$. If the linear independence constraint qualification (LICQ), the sufficient second order optimality conditions (SSOC) and the strict complementarity condition (SCC) are satisfied in this neighbourhood, then we have that
	\begin{itemize}
		\item $\controlstar(\cdot, \state, \disturbancezero)$ is an isolated local minimizer and the respective Lagrange multipliers are unique,
		\item for $(\statedisturbed, \disturbance )$ in a neighbourhood of $( \state, \disturbancezero )$ there exists a unique local minimizer $\controlstar(\cdot, \statedisturbed, \disturbance)$ which satisfies LICQ, SSOC and SCC and is differentiable with respect to $\statedisturbed$ and $\disturbance$,
		\item there exist a Lipschitz constant $L_\costfunction$ such that
		\begin{align}
			& | \costfunctionN( \statedisturbed, \feedbackdisturbedN(\statedisturbed, \disturbance), \disturbance) - \costfunctionN(\state, \feedbackN(\state, \disturbancezero), \disturbancezero) | \nonumber \\
			\label{Sensitivity theory:thm:sensitivity:eq1}
			& \qquad \leq L_\costfunction \left( \| \statedisturbed \|_{\state} + \| \disturbance \|_{\disturbancezero} \right)
		\end{align}
		holds and
		\item there exist a Lipschitz constant $L_\control$ such that for the updated control $\controldisturbed(\cdot, \statedisturbed, \disturbance)$ from \eqref{Sensitivity theory:eq:update formula}
		the following estimate holds:
		\begin{align*}
		& \distance_{\controlset} ( \controldisturbed(\statedisturbed, \disturbance), \controlstar(\state, \disturbancezero)) ) \leq L_\control \left( \| \statedisturbed \|_{\state} + \| \disturbance \|_{\disturbancezero} \right)
		\end{align*}
	\end{itemize}
\end{theorem}

Using this result, Assumption \ref{Stability:ass:sensitivity} can be verified.
\begin{proposition}\label{Sensitivity theory:prop:assumptions}
	If there exists a set $\updateset \subset \stateset$ containing $\statestar$ such that conditions of Theorem \ref{Sensitivity theory:thm:sensitivity} hold for all tupels $(\statedisturbed, \state, \disturbance, \disturbancezero)$ with $\state, \statedisturbed \in \updateset$, $\disturbancezero, \disturbance \in \disturbanceset$ satisfying Assumption \ref{Preliminaries:ass:bounded disturbance}, then Assumption \ref{Stability:ass:sensitivity} holds.
\end{proposition}
\textbf{Proof:} Follows directly from \eqref{Sensitivity theory:thm:sensitivity:eq1} and the fact that differentiability of $\stagecost$ implies the existence of a local Lipschitz constant $L_\stagecost$. \qed

Note that in order to apply both Proposition \ref{Sensitivity theory:prop:assumptions} and Theorem \ref{Stability:thm:performance} the set $\updateset$ is not necessarily large, a fact that otherwise may exclude such an approach in the presence of state and control constraints. 
In particular, Theorem \ref{Sensitivity theory:thm:sensitivity} requires the open--loop control structure to remain unchanged in a neighbourhood of the optimal solution despite changes in $\state$ and $\disturbancezero$. Most importantly, this property is only required locally. On a larger scale, the closed--loop control structure may change since we allow for an intermediate reoptimization. Hence, despite the fact that the update formula \eqref{Sensitivity theory:eq:update formula} is restricted to a certain neighbourhood of the open--loop solution, the MPC update approach using sensitivities is only locally restricted to that particular neighbourhood, i.e. for each visited closed--loop state this neighbourhood and the respective control structure may change.

%In the MPC context, we additionally don't require optimal controls in each step, see, e.g., Chapter 7.9 of \cite{GruenePannek2011} and references therein. Hence, considering a non--optimal control $\controlapprox(\cdot, \state, \disturbancezero)$ and defining $\valuefunctionNapprox(\state, \disturbancezero) := \costfunction(\state, \controlapprox, \disturbancezero)$ and $\feedbackNapprox(\state, \disturbancezero) := \controlapprox(0, \state, \disturbancezero)$, the descent property \eqref{Stability:thm:performance:eq1} can be relaxed to
%\begin{align*}
%	\valuefunctionNapprox(\state, \disturbancezero) \geq \valuefunctionNapprox(\dynamic(\state,\feedbackNapprox(\state, \disturbancezero), \disturbancezero), \disturbancezero) + \alpha \ell(\state, \feedbackNapprox(\state, \disturbancezero), \disturbancezero )
%\end{align*}
%
%Feasibility...
%
%\begin{theorem}\label{Sensitivity theory:thm:feasibility}
%	Was wenn nur fuer einen Schritt feasible + Abstieg
%\end{theorem}
%\textbf{Proof:} \qed

% ------------------------------------------------------------------------
\section{Numerical Results}
\label{Section:Numerical Results}
% ------------------------------------------------------------------------

\begin{figure}[!ht]
	\begin{center}
		\begin{tikzpicture}
			  %% front wheel
			  % mass top
			  \draw[fill] (4.5,4.5) circle (0.05) node[above right] {\relsize{-2}{$A$}};
			  \draw (4.25,4) -- (4.75,4);
			  \draw (4.5,4) -- (4.5,4.5);
			
			  % damper 
			  \draw (4.75,3) -- (4.75,3.4);
			  \draw (4.75,3.4) -- (4.6,3.4) -- (4.6,3.6);
			  \draw (4.75,3.4) -- (4.9,3.4) -- (4.9,3.6);
			  \draw (4.65,3.5) -- (4.85,3.5);
			  \draw (4.75,3.5) -- (4.75,4);
			  \draw (5.4,3.5) node {\relsize{-2}{$\control_1(t)$}};
			  
			  % spring
			  \draw (4.25,3) -- (4.25,3.25) -- (4.15,3.3) -- (4.35,3.4) -- (4.15,3.5) -- (4.35,3.6) -- (4.15,3.7) -- (4.25,3.75) -- (4.25,4);
			  \draw (3.85,3.5) node {\relsize{-2}{$k_3$}};
			  
			  % wheel
			  \draw (4.25,3) -- (4.75,3);
			  \draw (4.5,3) -- (4.5,1.5);
			  \draw (4.5,1.5) circle (1.1);
			  \draw[fill] (4.5,1.5) circle (0.05) node[above right] {\relsize{-2}{$m_1, \state_1$}};
			  \draw (4.25,1.5) -- (4.75,1.5);
			  
			  % damper
			  \draw (4.75,0.5) -- (4.75,0.9);
			  \draw (4.75,0.9) -- (4.6,0.9) -- (4.6,1.1);
			  \draw (4.75,0.9) -- (4.9,0.9) -- (4.9,1.1);
			  \draw (4.65,1.0) -- (4.85,1.0);
			  \draw (4.75,1.0) -- (4.75,1.5);
			  \draw (5.2,1.0) node {\relsize{-2}{$d_1$}};
			  
			  % spring
			  \draw (4.25,0.5) -- (4.25,0.75) -- (4.15,0.8) -- (4.35,0.9) -- (4.15,1) -- (4.35,1.1) -- (4.15,1.2) -- (4.25,1.25) -- (4.25,1.5);
			  \draw (3.85,1) node {\relsize{-2}{$k_1$}};
						  
			  \draw[->] (5,4.5) -- (5,4.1);
			  \draw (4.9,4.5) -- (5.1,4.5);
			  \draw (5.5,4.3) node {\relsize{-2}{$f_1$}};

			  %% back wheel
			  % mass top
			  \draw[fill] (0.5,4.5) circle (0.05) node[above right] {\relsize{-2}{$B$}};
			  \draw (0.25,4) -- (0.75,4);
			  \draw (0.5,4) -- (0.5,4.5);
			
			  % damper 
			  \draw (0.75,3) -- (0.75,3.4);
			  \draw (0.75,3.4) -- (0.6,3.4) -- (0.6,3.6);
			  \draw (0.75,3.4) -- (0.9,3.4) -- (0.9,3.6);
			  \draw (0.65,3.5) -- (0.85,3.5);
			  \draw (0.75,3.5) -- (0.75,4);
			  \draw (1.4,3.5) node {\relsize{-2}{$\control_2(t)$}};
			  
			  % spring
			  \draw (0.25,3) -- (0.25,3.25) -- (0.15,3.3) -- (0.35,3.4) -- (0.15,3.5) -- (0.35,3.6) -- (0.15,3.7) -- (0.25,3.75) -- (0.25,4);
			  \draw (-0.15,3.5) node {\relsize{-2}{$k_4$}};
			  
			  % wheel
			  \draw (0.25,3) -- (0.75,3);
			  \draw (0.5,3) -- (0.5,1.5);
			  \draw (0.5,1.5) circle (1.1);
			  \draw[fill] (0.5,1.5) circle (0.05) node[above right] {\relsize{-2}{$m_2, \state_2$}};
			  \draw (0.25,1.5) -- (0.75,1.5);
			  
			  % damper
			  \draw (0.75,0.5) -- (0.75,0.9);
			  \draw (0.75,0.9) -- (0.6,0.9) -- (0.6,1.1);
			  \draw (0.75,0.9) -- (0.9,0.9) -- (0.9,1.1);
			  \draw (0.65,1.0) -- (0.85,1.0);
			  \draw (0.75,1.0) -- (0.75,1.5);
			  \draw (1.2,1.0) node {\relsize{-2}{$d_2$}};
			  
			  % spring
			  \draw (0.25,0.5) -- (0.25,0.75) -- (0.15,0.8) -- (0.35,0.9) -- (0.15,1) -- (0.35,1.1) -- (0.15,1.2) -- (0.25,1.25) -- (0.25,1.5);
			  \draw (-0.15,1) node {\relsize{-2}{$k_2$}};
					  
			  % road
			  \draw (0.25,0.5) -- (0.75,0.5);
			  \draw (0.5,0.5) -- (0.5,0.4);
			  \draw (-1,0.5) sin (-0.5,0.6) cos (0,0.5) sin (0.5,0.4) cos (1,0.5) sin (1.5,0.6) cos (2,0.5) sin (2.5,0.4) cos (3,0.5) sin (3.5,0.6) cos (4,0.5) sin (4.5,0.4) cos (5,0.5) sin (5.5,0.6) cos (6,0.5);
			  \draw (0.5,0.0) node {\relsize{-2}{$\disturbancezero(t - \Delta)$}};
			  \draw (4.25,0.5) -- (4.75,0.5);
			  \draw (4.5,0.5) -- (4.5,0.4);
			  \draw (4.5,0.0) node {\relsize{-2}{$\disturbancezero(t)$}};
			  
			  \draw[->] (-0.0,4.5) -- (-0.0,4.1);
			  \draw (-0.1,4.5) -- (0.1,4.5);
			  \draw (-0.5,4.3) node {\relsize{-2}{$f_2$}};

			 % car
			\draw (0.5, 4.5) -- (4.5, 4.5);
			\draw[fill] (2.5,4.5) circle (0.05) node[above] {\relsize{-2}{$m_3, I, \state_3$}};
			\draw (3.5,4.3) node {\relsize{-2}{$a$}};
			\draw (1.5,4.7) node {\relsize{-2}{$b$}};
			\draw (1.0, 4.2) -- (4.0, 4.8);
			\draw (4.0,4.65) node {\relsize{-2}{$\state_4$}};
		\end{tikzpicture}
		\caption{Schematical sketch of a halfcar subject to road excitation $\disturbancezero$}
		\label{Numerical Results:fig:halfcar}
	\end{center}
\end{figure}
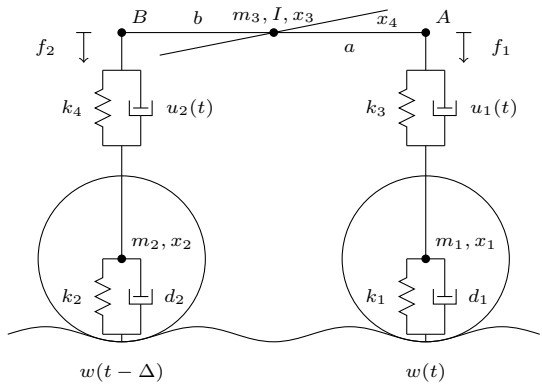

To illustrate our results, we consider a halfcar model given from \cite{SpeckertDresslerRuf2009, PoppSchiehlen2010} with proactive dampers given by the second order dynamics
\begin{align*}
	m_1 \ddot{\state}_1 & = m_1 g + f_3 - f_1 \\
	m_2 \ddot{\state}_2 & = m_2 g + f_4 - f_2 \\
	m_3 \ddot{\state}_3 & = m_3 g - f_3 - f_4 \\
	I \ddot{\state}_4 & = \cos(\state_4) ( b f_3 - a f_4 )
\end{align*}
where the control enters the forces
\begin{align*}
	f_i & = k_i ( \state_i - \disturbancezero_i ) + d_i (\dot{\state}_i - \dot{\disturbancezero}_i), \quad i = 1, 2 \\
	f_3 & = k_3 ( \state_3 - \state_1 - b \sin(\state_4) ) + \control_1 ( \dot{\state}_3 - \dot{\state}_1 - b \dot{\state}_4 \cos(\state_4) ) \\
	f_4 & = k_4 ( \state_3 - \state_2 + a \sin(\state_4) ) + \control_2 ( \dot{\state}_3 - \dot{\state}_2 + a \dot{\state}_4 \cos(\state_4) )
\end{align*}
see Fig. \ref{Numerical Results:fig:halfcar} for a schematical sketch.
Here, $\state_1$ and $\state_2$ denote the centers of gravity of the wheels, $\state_3$ the respective center of the chassis and $\state_4$ the pitch angle of the car. The disturbances $\disturbancezero_1$, $\disturbancezero_2$ are connected via $\disturbancezero_1(t) = \disturbancezero(t)$, $\disturbancezero_2(t) = \disturbancezero(t - \Delta)$. The remaining constants of the halfcar are displayed in Tab. \ref{Numerical Results:tab:halfcar}.
\begin{table}[!ht]
	\begin{center}
		\begin{tabular}{|c|c|c|c|} \hline
			name & symbol & quantity & unit\\ \hline\hline
			distance to joint & $a, b$ & $1$ & $m$ \\
			mass wheel & $m_1, m_2$ & $15$ & $kg$ \\
			mass chassis & $m_3$ & $750$ & $kg$\\
			inertia & $I$ & $500$ & $kg\, m^2$\\
			spring constant wheels & $k_1, k_2$ & $2\cdot 10^5$ & $kN/m$\\
			damper constant wheels & $d_1, d_2$ & $2\cdot 10^2$ & $kNs/m$\\
			spring constant chassis & $k_3, k_4$ & $1 \cdot 10^5$ & $kN/m$\\
			gravitational constant & $g$ & $9.81$ & $m/s^2$ \\
		\hline
		\end{tabular}
		\caption{Parameters for the halfcar example}
		\label{Numerical Results:tab:halfcar}
	\end{center}
\end{table}
For this problem, we apply MPC cost functional
\begin{align*}
	\costfunctionN(\state, \control, \disturbancezero) := \sum_{\openloopindex = 0}^{N-1} \mu_R F_R( \openloopindex ) + \mu_A F_A( \openloopindex )
\end{align*}
following ISO 2631 with horizon length $N = 5$. The handling objective is implemented via
\begin{small}
\begin{align*}
	& F_R( \openloopindex ) := \int\limits_{\openloopindex T}^{(\openloopindex+1) T} \left( \frac{[ k_1 ( \state_1(t) - \disturbancezero_1(t) ) + d_1 ( \dot{\state}_1(t) - \dot{\disturbancezero}_1(t) ) ] - F_1}{F_1}  \right)^2 \\
			& + \left( \frac{[ k_2 ( \state_2(t) - \disturbancezero_2(t) ) + d_2 ( \dot{\state}_2(t) - \dot{\disturbancezero}_2(t) ) ] - F_2}{F_2}  \right)^2\, dt
\end{align*}
\end{small}
with nominal forces $F_1 = ( a \cdot g \cdot ( m_1 + m_2 + m_3 ) ) / ( a + b )$, $F_2 = ( b \cdot g \cdot ( m_1 + m_2 + m_3 ) ) / ( a + b )$ whereas minimizing the chassis jerk
\begin{small}
\begin{align*}
	F_A(\openloopindex) := \int\limits_{\openloopindex T}^{(\openloopindex+1)T} \left( m_3 \dddot{\state}_3(t) \right)^2 \, dt
\end{align*}
\end{small}%
is used to treat the comfort objective. Both integrals are evaluated using a constant sampling rate of $T = 0.1s$ during which the control are held constant, i.e. the control is implmented in a zero--order hold manner. Additionally, the control constraints $\controlsetconstrained = [ 0.2 kNs/m, 5 kNs/m ]^2$ limit the range of the controllable dampers. Within the MPC scheme outlined in Algorithm \ref{Sensitivity theory:alg:mpc}, we compute the nominal disturbance $\disturbancezero(\cdot)$ and the corresponding derivates from road profile measurements taken at a sampling rate of $0.002s$ via a fast Fourier transformation (FFT). Both the states of the system and the road profile measurements are modified using a disturbance which is uniformly distributed in the interval $[-0.025, 0.025]$.
\begin{figure}[!ht]
	\includegraphics[width=0.48\textwidth]{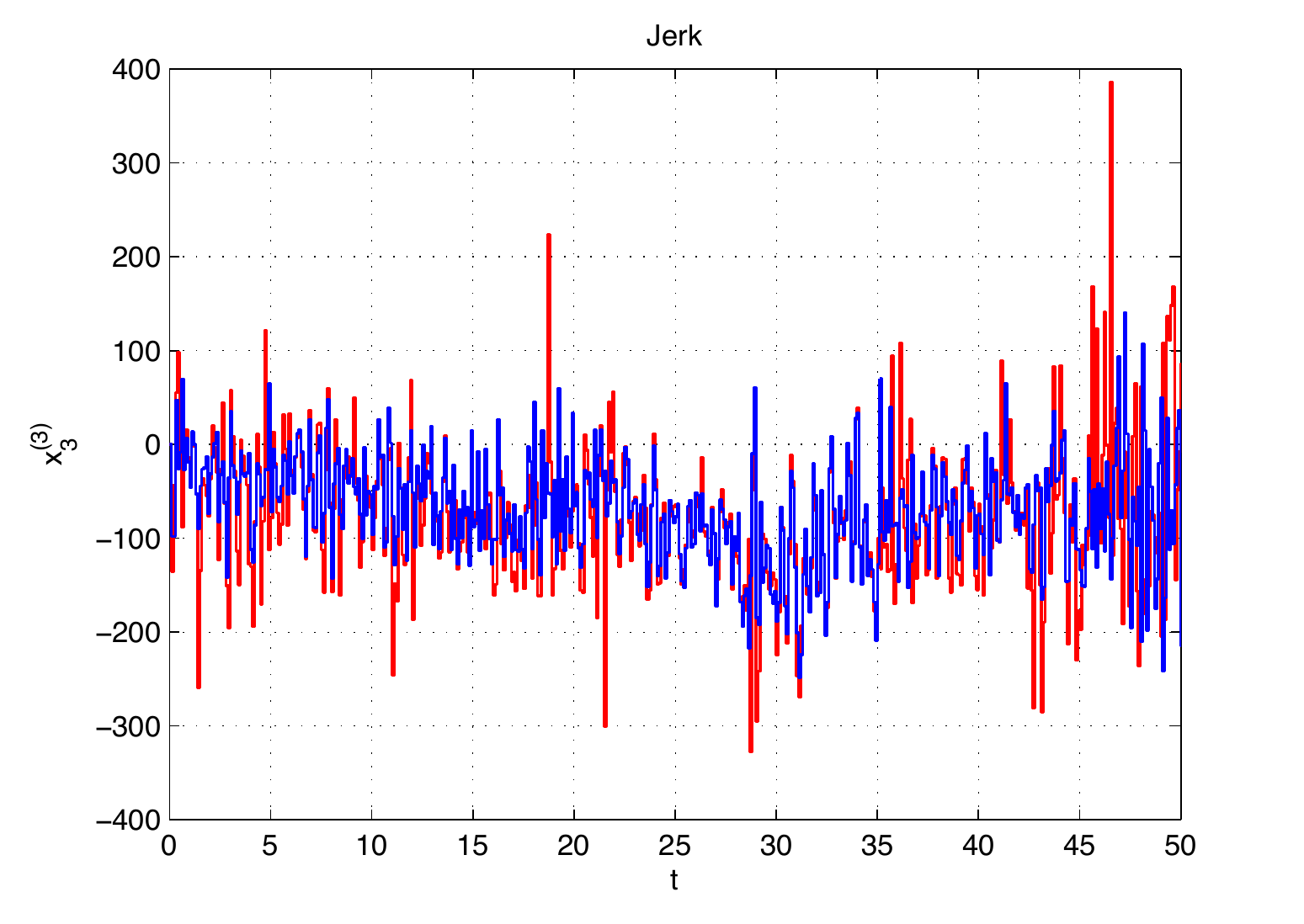}
	\caption{Chassis jerk for MPC with (blue) and without sensitivity update (red)}
	\label{Numerical Results:fig:jerk}
\end{figure}
%\begin{figure}[!ht]
%	\includegraphics[width=0.48\textwidth]{handling}
%	\caption{Handling for MPC with (blue) and without sensitivity update (red)}
%	\label{Numerical Results:fig:handling}
%\end{figure}

As expected, the updated control law shows a better performance which can not only be observed from Fig. \ref{Numerical Results:fig:jerk}% and \ref{Numerical Results:fig:handling}
, but also in terms of the closed loop costs: For the chosen industrial road data we obtain an improvement of approximately $8.2 \%$ using the sensitivity update \eqref{Sensitivity theory:eq:update formula}. Although this seems to be a fairly small improvement, the best possible result obtained by a full reoptimization reveals a reduction of approximately $10.5 \%$ of the closed loop costs.

We like to note that due to the presence of constraints it is a priori whether the conditions of Theorem \ref{Sensitivity theory:thm:sensitivity} hold at each visited point along the closed loop. Such an occurrance can be detected online by checking for violations of constraints or changes in the control structure. Yet, due to the structure of the MPC algorithm, such an event has to be treated if one of the constraints is violated at open loop time instant $\openloopindex = 1$ only which was not the case for our example.

% ------------------------------------------------------------------------
\section{Conclusions and Outlook}
\label{Section:Conclusions and Outlook}
% ------------------------------------------------------------------------

We considered MPC without stabilizing terminal constraints or Lyapunov type terminal costs in the presence of disturbances. For this setting, we presented stability and performance results for abstract control updates. The required assumptions are shown to be fulfilled if sensitivity updates on the control can be performed.

Future research concerns the analysis of other types of control updates such as realtime iterations which have been considered for MPC with stabilizing terminal constraints and costs. Using the presented result, we hope to obtain a unified stability and performance analysis of MPC in the presence of disturbances and respective control updates. Regarding the halfcar example, the influence of other types of interpolation/extrapolation methods on the MPC performance such as road profiles based on statistical data will be analysed.

% ------------------------------------------------------------------------
\begin{ack}                              
	This work was partially funded by the German Federal Ministry of Education and Research (BMBF), grant no. 03MS633G.
\end{ack}
% ------------------------------------------------------------------------

% ------------------------------------------------------------------------
\bibliographystyle{alpha}        % Include this if you use bibtex 
\bibliography{ifacconf}
% ------------------------------------------------------------------------

%% ------------------------------------------------------------------------
%\appendix
%% ------------------------------------------------------------------------

% ------------------------------------------------------------------------
\end{document}